\documentclass[12pt,a4paper]{amsart}

\usepackage{color}
\input{xypic}

\usepackage[a4paper]{geometry}
\usepackage[all]{xy}
\usepackage{amsthm, amssymb, amsmath}
\usepackage{amsopn}

\usepackage{mathrsfs}
\usepackage{enumitem}
\usepackage{mathtools}
\usepackage{tikz-cd}

\usepackage{hyperref}
\hypersetup{
    colorlinks=true,
    linkcolor=blue,
    filecolor=magenta,
    urlcolor=cyan,
    }

\advance \topmargin by -\headheight
\advance \topmargin by -\headsep
\evensidemargin \oddsidemargin
\numberwithin{equation}{section}
\setenumerate[1]{label={\alph*)}} % Global setting
\newcommand{\al}{\alpha}

\newcommand{\D}{\Delta}

\newcommand{\La}{\Lambda}
\newcommand{\la}{\lambda}

\def\C{\mathbb{C}}

\newcommand{\Q}{\mathbb{Q}}
\newcommand{\Z}{\mathbb{Z}}
\newcommand{\mbo}{\mathbb{O}}

\newcommand{\CL}{\mathcal{L}}

\newcommand{\g}{\mathfrak{g}}
\newcommand{\h}{\mathfrak{h}}
\newcommand{\sll}{\mathfrak{sl}}

\newcommand{\ov}[1]{\overline{#1}}
\newcommand{\what}[1]{\widehat{#1}}
\newcommand{\wtil}[1]{\widetilde{#1}}

\newcommand{\vac}{\mathbf{1}}
\newcommand{\vir}{\text{Vir}}

\DeclareMathOperator{\ad}{ad}

\DeclareMathOperator{\im}{im}

\swapnumbers
\newtheoremstyle{exps}{\topsep}{\topsep}{}{0pt}{\bfseries}{.}{0pt}{}

\DeclareMathOperator{\heit}{ht}

\newtheorem*{thm*}{Theorem}
\newtheorem*{prop*}{Proposition}
\newtheorem*{lem*}{Lemma}
\newtheorem*{cor*}{Corollary}
\newtheorem*{rem*}{Remark}
\newtheorem{thm}{Theorem}[section]

\theoremstyle{definition}
\newtheorem*{defn*}{Definition}
\newtheorem*{exer*}{Exercise}

\newtheorem*{problem*}{Problem}
\newtheorem{rem}[thm]{Remark}

\theoremstyle{exps}

\title[]{Some new isomorphisms of affine $W$-algebras}

\begin{document}

\begin{center}
{\LARGE \bf Some isomorphisms between exceptional $W$-algebras} \par \bigskip

\renewcommand*{\thefootnote}{\fnsymbol{footnote}}
{\normalsize Jethro van Ekeren\footnote{email: \texttt{jethro@impa.br}}} \\
{\footnotesize \emph{Instituto de Matem\'{a}tica Pura e Aplicada, Rio de Janeiro, RJ 22.460-320, Brazil}} \\
\vspace{2mm}
{\normalsize Shigenori Nakatsuka}\footnote{email: \texttt{shigenori.nakatsuka.2022@gmail.com}} \\
{\footnotesize \emph{Research Institute for Mathematical Sciences, Kyoto University, Kyoto 606-8502, Japan}}

\par \bigskip
%\par
%
\par \bigskip

\end{center}

%\vspace*{10mm}

\noindent
\textbf{Abstract.} We prove some isomorphisms between exceptional $W$-algebras associated with exceptional simple Lie algebras.

\vspace*{5mm}

%\tableofcontents

\section{Introduction}

The purpose of this note is to prove three isomorphisms of vertex algebras, namely those listed in Theorem \ref{thm:isoms.intro} below. All the vertex algebras in question are examples of $W$-algebras, by which we mean vertex algebras obtained via quantised Drinfeld-Sokolov reduction of affine vertex algebras \cite{FF} \cite{KRW}. These form a large class of vertex algebras with important applications in representation theory and theoretical physics.

The universal affine vertex algebra associated with a simple Lie algebra $\g$ and level parameter $k$ is denoted $V^k(\g)$. The reduction takes as an additional parameter a nilpotent element $f \in \g$ and the resulting vertex algebra $W^k(\g, f)$, known as the universal $W$-algebra, depends up to isomorphism on $f$ only through its adjoint orbit $\mbo = G \cdot f$.

For generic $k$ the universal $W$-algebra $W^k(\g, f)$ is simple, while at special values of $k$ the structure of its simple quotient can be very difficult to ascertain. For special levels of the form $k = -h^\vee + p/u$, introduced and studied by Kac and Wakimoto and called admissible levels \cite{KW88}, the simple quotient $L_k(\g)$ of $V^k(\g)$ is in a certain sense ``partially integrable'' as a $\what\g$-module. This partial integrability allows much to be said about the structure and representation theory of the simple quotient of $W^k(\g, f)$, which we denote by $W(\mbo, p/u)$.

We prove the following theorem.
\begin{thm}\label{thm:isoms.intro}
There exist isomorphisms of vertex algebras
\begin{align}\label{eq:isoms.intro}
\begin{split}
\textup{(1)}\quad W(E_7(a_1), 19/16) &\cong W(B_3, 8/7), \\
\textup{(2)}\quad W(E_7(a_1), 18/17) &\cong W(G_2, 4/17), \\
\textup{(3)}\quad W(E_8(a_1), 31/27) &\cong W(G_2, 9/7).
%\text{(4)}\quad W(E_8(a_1), 31/28) &\cong W(F_4, 14/13) \\
%\text{(5)}\quad W(E_8(a_1), 30/29) &\cong W(E_8, 29/24)
\end{split}
\end{align}
\end{thm}
In the theorem statement we have used the Bala-Carter labeling of nilpotent orbits, so $E_7(a_1)$ is the subregular nilpotent orbit in $E_7$, etc., and on the right hand side the orbits are all principal. We expect a fourth isomorphism
\begin{align}\label{eq:isom.F4}
\text{(4)}\quad W(E_8(a_1), 31/28) &\cong W(F_4, 14/13),
%\text{(5)}\quad W(E_8(a_1), 30/29) &\cong W(E_8, 29/24)
\end{align}
whose proof should follow by the method employed in this article. The only obstruction is that certain OPE computations are too burdensome to carry out at the present time.

The $W$-algebras listed above are all examples of rational lisse vertex algebras.  For each integer $u \geq 1$ there is a nilpotent orbit $\mbo_u \subset \g$ for which $W(\mbo_u, p/u)$ is lisse \cite{Arakawa.assoc.var}. In reference to the notion of exceptional pair introduced in \cite{KW08.Rationality}, we refer to such simple $W$-algebras as exceptional.

In \cite{AE.JEMS} T. Arakawa and the first author have established rationality of a large class of exceptional $W$-algebras (see \cite{McRae} in general), and have computed their fusion rules in some cases. E. Rowell \cite{Rowell.private} noticed relations between these fusion rules and those of certain simple affine vertex algebras at positive integer level: $W(E_8(a_1), 31/27)$ is related to $L_5(G_2)$ for instance. These coincidences are explained by combining the isomorphisms listed above with a calculation of fusion rules of principal $W$-algebras, due to Frenkel, Kac and Wakimoto \cite{FKW}. Roughly speaking, the Grothendieck ring of the principal $W$-algebra $W(\g, p/u)$ factorises into a product of the Grothendieck rings of $L_{p-h^\vee}(\g)$ and $L_{u-h^\vee}(\g)$. Therefore taking either $p$ or $u$ close to $h^\vee$ yields a principal $W$-algebra whose fusion rules are similar to those of a simple affine vertex algebra. With this picture in mind, confirmed by computation of central charges and asymptotic growths, the isomorphisms of Theorem \ref{thm:isoms.intro} were discovered. They were then proved by developing character formulas, and using these together with facts about the representation theory of the Virasoro algebra, and explicit OPE computations.

Isomorphism of $W$-algebras has been investigated from several perspectives, and using a variety of methods. The Feigin-Frenkel isomorphism \cite{FF} between Langlands dual principal $W$-algebras is a very important example. Isomorphisms related to the triality phenomenon in theoretical physics \cite{GG} \cite{GR} have been studied and proved in \cite{Linshaw-W-infty}, \cite{CL-trialities} and \cite{CL-ortho-trialities}. The notion of collapsing level has been introduced and studied in \cite{collapsing-I} and \cite{collapsing-II}, motivated by applications \cite{collapsing-applic} to the representation theory of $\what \g$. Classification of admissible collapsing levels has been approached in \cite{AEM-collapsing} using associated varieties and asymptotics of characters, and using similar techniques, some other isomorphisms between classes of $W$-algebras were also proved. See also \cite{AE.JEMS} and \cite{Fasquel}. Characters are also fundamental to the approach used in this note, not through their asymptotic properties, but rather to constrain OPE structures in low conformal weight.

\emph{Acknowledgements:} These isomorphisms were discovered and proved in October 2023 while the authors attended ``Geometric Representation Theory and $W$-algebras'' in Edinburgh, and the results have been presented in December 2023 at the Indam workshop ``Vertex algebras, infinite dimensional Lie algebras and related topics'' in Rome. We are grateful to the organisers of these events. We would also like to thank T. Arakawa, A. Linshaw and K. Sun for useful suggestions and discussions. JvE has been supported by grant numbers Serrapilheira -- 2023-0001, CNPq 306498/2023-5 and FAPERJ 201.445/2021. SN has been supported by JSPS Overseas Research Fellowships Grant Number 202260077.

\section{Background}

\subsection{Vertex algebras}{\ }\label{subsec:va}

We recall that a vertex algebra (see \cite{va.beginners}, \cite{DSK06}) consists of a vector space $V$, a distinguished vector $\vac$ called the vacuum vector, an endomorphism $T : V \rightarrow V$ called the translation operator, and a collection of bilinear products $(a, b) \mapsto a_{(n)}b$ for $n \in \Z$. The formal series
\[
Y(a, z) = \sum_{n \in \Z} a_{(n)} z^{-n-1}
\]
is called the field associated with the state $a \in V$ (the coefficients $a_{(n)}$ are endomorphisms of $V$). These data are required to satisfy (1) $Y(a, z)b \in V((z))$ for each pair $a, b \in V$, (2) $Y(\vac, z) = I_V$, (3) $T\vac = 0$ and $[T, Y(a, z)] = \partial_z Y(a, z)$ for all $a \in V$, (4) $Y(a, z)\vac \in V[[z]]$ and $Y(a, z)\vac|_{z=0} = a$ for all $a \in V$, and (5) for each pair $a, b \in V$ there exists $N$ such that $(z-w)^N [Y(a, z), Y(b, w)] = 0$.

The vertex algebras we shall be concerned with in this paper are conformal. A conformal structure on a vertex algebra $V$ is a grading $V = \bigoplus_{n \in \Z_+} V_n$ by finite dimensional subspaces $V_n$, and a vector $L \in V_2$ such that if $Y(L, z) = \sum_{n \in \Z} L_n z^{-n-2}$ then (1) the coefficients $L_n$ furnish $V$ with a representation of the Virasoro Lie algebra
\begin{align*}
\CL = \bigoplus_{n \in \Z} \C L_n \oplus \C C, \quad [L_m, L_n] = (m-n) L_{m+n} +\delta_{m, -n} \frac{m^3-m}{12} C, \quad [C, \CL] = 0,
\end{align*}
in which $C$ acts as $c I_V$ for a constant $c$ called the central charge, (2) $T = L_{-1}$, and (3) the subspace $V_n$ is the eigenspace of $L_0$ of eigenvalue $n$. If $a \in V_\Delta$ then it is a straightforward consequence of the definitions that $a_{(n)}V_k \subset V_{k+\Delta-n-1}$. We call $\Delta$ the conformal weight of $a$.

It is common to denote $a_{(-1)}b$ as $:\!ab\!:$, and refer to this term as the normally ordered product of $a$ and $b$. The coefficients $a_{(n)}b$ for $n \geq 0$ form what is called the operator product expansion (OPE) of $a$ and $b$. It is convenient to gather these terms into the $\lambda$-bracket
\[
[a_\la b] = \sum_{n \in \Z_+} \frac{\la^n}{n!} a_{(n)}b.
\]
In any vertex algebra $V$ the following commutator formula holds for all $a, b, c \in V$:
\begin{align}\label{eq:comm.fla}
[a_{\la}[b_{\mu} c]] - [b_{\mu}[a_{\la}c]] = [[a_\la b]_{\la+\mu} c],
\end{align}
as well as further relations involving $\lambda$-brackets and normally ordered products, such as the noncommutative Wick formula {\cite[Section 3.3]{va.beginners}}
\begin{align}\label{eq:noncomm.Wick}
[a_{\la}{:bc:}] = :[a_\la b]c: + :b[a_\la c]: + \int_{0}^\la [[a_\la b]_\mu c] \, d\mu.
\end{align}

We say $\{a^1, \ldots, a^s\} \subset V$ is a strong generating set if the terms
\begin{align}\label{eq:str.gen.monom}
a^{i_1}_{(-n_1-1)} a^{i_2}_{(-n_2-1)}  \cdots a^{i_r}_{(-n_r-1)}  \vac, \quad \text{$r \geq 0$ and $n_i \geq 0$}
\end{align}
span $V$. From the axioms involving $\vac$ and $T$ it follows that $(Ta)_{(n)} = -n a_{(n-1)}$ for all $a \in V$ and all $n \in \Z$, so one could equally say that $V$ is strongly generated by a set $S$ if it is spanned by normally ordered polynomials in arbitrary derivatives of elements of $S$.

In Section \ref{sec:isoms} below we shall use \eqref{eq:comm.fla} to constrain possible OPE relations in vertex algebras in which we have good control of the conformal weights of generators and relations.

%%%%%%%%%%%%%%%%%%%%%%%%%%%%%

\subsection{Representation theory of the Virasoro Lie algebra}{\ }\label{sec: Representation theory of the Virasoro Lie algebra}

We review a few standard facts about highest weight modules over the Virasoro Lie algebra $\CL$ that we will use in Section \ref{sec:isoms}. See \cite{Feigin-Fuchs}, \cite{Iohara-Koga} or \cite{BNW}. We denote by $M(c, h)$ the Verma module $U(\CL) \otimes_{U(\CL_+)} \C v$ in which $\CL_+ = \bigoplus_{n \in \Z_{\geq 0}} \C L_n \oplus \C C$ acts on $v$ as follows: $C v = cv$, $L_0v = hv$ and $L_{>0} v = 0$. We denote by $L(c, h)$ the irreducible quotient of $M(c, h)$, and for $h = 0$ we have the the vacuum module $V(c) = M(c, h) / U(\CL) L_{-1}v$.

Let
\begin{align*}
\nu = \frac{1}{12}\left(c - 13 + \sqrt{(c-1)(c-25)} \right) \qquad \beta = \sqrt{(\nu+1)^2 - 4h \nu}.
\end{align*}
In general the equation
\[
r + \nu s + \beta = 0
\]
might have zero, one or infinitely many integer solutions $(r, s)$. For all the values of $c$ relevant to this paper, we have $\nu \notin \Q$ and therefore in each such case only zero or one solution $(r, s)$. In such cases the only Verma module linked with $M(c, h)$ (besides itself) is $M(c, h+rs)$.

\subsection{Simple Lie algebras and nilpotent orbits}{\ }

Let $\g$ be a finite dimensional simple Lie algebra of rank $\ell$, with Cartan subalgebra $\h$ and root system $\ov{\D} \subset \h^*$. A choice of triangular decomposition yields a set $\ov\D_+$ of positive roots, and its subset $\Pi = \{\ov\al_1, \ldots, \ov\al_\ell\}$ of simple roots. The simple coroots are $\ov\al_i^\vee = 2 \nu^{-1}(\ov\al_i) / (\ov\al_i, \ov\al_i)$, where $\nu : \h \rightarrow \h^*$ is the linear isomorphism induced by the choice of non degenerate invariant bilinear form $(\cdot, \cdot) : \g \times \g \rightarrow \C$ which, by convention, is normalised so that the highest root $\theta$ has squared norm $2$. The fundamental weights $\{\varpi_i\}$ are the elements of the basis of $\h^*$ canonically dual to the basis $\{\ov\al_i^\vee\}$ of $\h$. It is well known that the Weyl vector, defined as $\ov\rho = \frac{1}{2}\sum_{\ov\al \in \ov\D_+} \ov\al$, coincides with the sum $\sum_{i=1}^\ell \varpi_i$.

Now let $f \in \g$ be an $\ad$-nilpotent element. We choose an $\sll_2$-triple $\{e, h, f\}$ containing $f$ (the existence of which is guaranteed by the Jacobson-Morozov theorem), and we consider the associated Dynkin grading $\g = \bigoplus_{j \in (1/2)\Z} \g_j$ of $\g$ by eigenspaces of $\ad(h)$, namely $[h, x] = 2jx$ for $x \in \g_j$. A nilpotent orbit is said to be even if the Dynkin grading of any of its elements is integral. In this note we will be concerned with a few specific nilpotent orbits, all of which are even. The centraliser $\g^f = \{x \in \g \mid [f, x] = 0\}$ is a graded subalgebra of $\g$ with respect to the Dynkin grading. By $\sll_2$-theory $\dim(\g_{-j}) = \dim(\g_j)$, the graded components of $\g^f$ are concentrated in non-positive degrees, and $\dim(\g^f_{-j}) = \dim(\g^e_j) = \dim(\g_{j}) - \dim(\g_{j+1})$ for $j \geq 0$. So it is possible to read off the dimensions $\dim(\g_j)$ from the list of degrees of a graded basis of $\g^e$. In the principal case each degree appears at most once, and the degrees are precisely the exponents of the Lie algebra. We record the degrees for some relevant principal orbits
\begin{align}\label{eq:gr.poly.prin}
\begin{split}
{G_2} : {} & (1, 5), \\
{B_3} : {} & (1, 3, 5), \\
{F_4} : {} & (1, 5, 7, 11),
\end{split}
\end{align}
and subregular orbits
\begin{align}\label{eq:gr.poly.subreg}
\begin{split}
{E_6(a_1)} : {} & (1, 2, 3, 4, 5, 5, 7), \\
{E_7(a_1)} : {} & (1, 3, 5, 5, 7, 8, 9, 11, 13), \\
{E_8(a_1)} : {} & (1, 5, 7, 9, 11, 13, 14, 17, 19, 23).
\end{split}
\end{align}
Here and throughout the paper we denote a nilpotent orbit by its Bala-Carter label \cite{Bala-Carter-2} {\cite[Chapter 8]{Coll-McG}}. Since we may take a Cartan subalgebra of $\g$ which contains $h$, it is evidently possible to arrange that the root space decomposition of $\g$ be compatible with a given Dynkin grading. We fix such a choice and write $\ov\D_{>0}$ for the set of roots associated with $\g_{>0} = \bigoplus_{j > 0} \g_j$.

\subsection{Admissible weights and the Kac-Wakimoto character formula}{\ }\label{subsec:KW-char-fla}

For $\g$ as above we consider the untwisted affine Kac-Moody algebra $\what{\g} = \g[t, t^{-1}] \oplus \C K \oplus \C d$, where
\[
[at^m, bt^n] = [a, b]t^{m+n} + m \delta_{m, -n} (a, b) K, \quad [K, \what\g] = 0, \quad [d, at^m] = m at^m,
\]
with Cartan subalgebra $\what\h = \h \oplus \C K \oplus \C d$. We write weights $\lambda \in \what\h^*$ relative to a decomposition $\what\h^* = \h^* \oplus \C \delta \oplus \C \Lambda_0$ where $\La_0(\h + \C d) = \delta(\h + \C K) = 0$ and $\La_0(K) = \delta(d) = 1$ (see \cite{IDLA}).

The corresponding root system is the union of the sets $\what{\D}_+^{\text{re}} = \{\ov\al + n \delta \mid n \in \Z\}$ of real roots, and $\{n \delta \mid n \in \Z \backslash \{0\}\}$ of imaginary roots. Real coroots are defined as in the finite type case, and the corresponding fundamental weights are $\{\La_0, \La_1, \ldots, \La_\ell\}$ where $\La_i = \varpi_i + a_i^\vee \La_0$. Here the integers $a_i^\vee = \left<\varpi_i, \theta^\vee\right>$ are the comarks, which can be read from {\cite[pp. 54-55, p. 79]{IDLA}}.

A weight $\la \in \what\h^*$ is said to be dominant integral if
\begin{align*}
\left<\la, \al_i^\vee\right> \in \Z_{\geq 0}, \quad \text{$i = 0, 1, \ldots, \ell$},
\end{align*}
and the set $P_+$ of such weights is just $\Z_+ \{\La_0, \ldots, \La_\ell\}$. The level of a weight $\la$ is the number $\left<\la, K\right>$. For example the affine Weyl vector $\rho = \sum_{i=0}^\ell \La_i$ has level $h^\vee = 1 + \sum_{i=1}^\ell a_i^\vee$ the dual Coxeter number of $\g$. Indeed $\rho = h^\vee \La_0 + \ov\rho$.

The formal character of the irreducible $\what\g$-module of dominant integral highest weight $\la$ is given by the Weyl-Kac character formula
\begin{align*}
\chi_{L(\la)} = \frac{A_{\la+\rho}}{D}, \quad \text{where} \quad A_\la = \sum_{w \in \what{W}} \varepsilon(w)  e^{w(\la)} \quad \text{and} \quad D = e^{\rho} \prod_{\alpha \in \what\D_+} (1-e^{-\alpha}).
\end{align*}
The affine Weyl group $\what{W}$, which appears in these formulas, is the semidirect product $W \ltimes t_{Q^\vee}$, where $W$ is the finite Weyl group and $t_{\al}$ is the translation acting on $\what\h^*$ via
\begin{align*}
t_\al(\la) = \la + \lambda(K)\al - \left[ (\al, \la) + \frac{|\al|^2}{2} \la(K) \right] \delta
\end{align*}
In this context one has the Weyl denominator formula
\begin{align}\label{eq:denom.fla}
A_\rho = D.
\end{align}
The extended affine Weyl group $\wtil W$, which will also appear below, is $W \ltimes t_{P^\vee}$ where $P^\vee = Q^*$. The quotient $\wtil W / \what W$ is a finite abelian group isomorphic to $P/Q^\vee$.

The class of admissible weights $\la \in \what\h^*$ was introduced by Kac and Wakimoto in \cite{KW88}, and it was proven in that article that $\chi_{L(\la)}$ admits a formula similar to the Weyl-Kac character formula for $\la$ admissible.

For a weight $\la \in \what{\h}^*$ we consider the subset of coroots
\[
\D_+^{\vee, \text{re}}(\la) = \{\al^\vee \in \D_+^{\vee, \text{re}} \mid \la(\al^\vee) \in \Z \}
\]
and $\Pi^\vee(\la)$ the system of simple coroots relative to $\D_+^{\vee, \text{re}}(\la)$. The weight $\la$ is said to be admissible if it is dominant integral with respect to $\Pi^\vee(\la)$. The level $k$ is said to be principal admissible if $k\La_0$ is admissible and furthermore
\begin{align*}
\Pi^\vee(k\La_0) = \{\al_1^\vee, \ldots, \al_\ell^\vee, -\theta^\vee + uK\}.
\end{align*}
(If $\g$ is simply laced then this is, in fact, the only possible form that $\Pi^\vee(k\La_0)$ can take. More generally $\Pi^\vee(k\La_0)$ might take a similar form in which $\theta^\vee$ is replaced by $\theta^\vee_{\text{long}}$, in which case $k$ is said to be coprincipal admissible.) Concretely the principal admissible levels are $k = -h^\vee + p/u$ where $p, u$ are coprime positive integers with $p \geq h^\vee$ and $u$ coprime to the lacety $r^\vee$ of $\g$. We say that a weight $\la$ is principal admissible of level $k$ if $\Pi^\vee(\la) = y\left(\Pi^\vee(k\La_0)\right)$ for some $y \in \wtil{W}$. The principal admissible weights corresponding to a given choice $y = t_{\beta} \ov{y}$, where $\ov y \in W$ and $\beta \in {P^\vee}$, are of the form
\begin{align}\label{eq:adm.param}
\la = y(\phi(\nu)) - \rho, \quad \nu \in P_+^{p, \text{reg}}.
\end{align}
Here $\phi : \what\h^* \rightarrow \what\h^*$ is the isometry $\phi(k\La_0 + \ov\la + t\delta) = (k/u)\La_0 + \ov\la + ut\delta$ and {{$P_+^{p, \text{reg}}$ is the set of dominant integral weights at level $p$}}.

Upon coordinatising $\what\h^*$ as follows
\begin{align*}
(\tau, z, t) = -2\pi i \left( \tau \La_0 + \sum_{i=1}^\ell z_i \al_i + t\delta \right),
\end{align*}
Kac and Wakimoto obtain the following formula for the character of $L(\la)$, for $\la$ principal admissible as above:
\begin{align}\label{eq:adm.char.fla}
\chi_{L(\la)}(\tau, z, t) = \frac{A_{\nu}(u \ov y^{-1} t_{-\beta/u}(\tau, z/u, t/u^2))}{A_\rho(\tau, z, t)}.
\end{align}
In fact we shall typically suppress the coordinate $t$, and write
\begin{align*}
q = e^{-2\pi i \tau \delta}, \quad y_i = e^{2\pi i z_i \ov\al_i} \quad \text{and} \quad \text{$y_\al = \prod_{i} y_i^{k_i}$ for $\al = \sum_i k_i \ov\al_i \in \ov\D_+$},
\end{align*}
so in particular the Weyl denominator is
\begin{align}\label{eq:denom.qy}
D(q, y_i) = \prod_{n=1}^\infty (1-q^n)^\ell \prod_{\al \in \ov\D_{+}} (1-y_\alpha^{-1} q^{n-1}) (1-y_\alpha q^{n}).
\end{align}

To illustrate we consider $\g = \sll_2$ and $k = -2 + 2/5$, so the following two weights are principal admissible of level $k$:
\begin{align*}
\la_0 = t_{-\varpi}(\phi(\rho)) - \rho \quad \text{and} \quad \la_1 = t_{-2\varpi}(\phi(\rho)) - \rho.
\end{align*}
The Weyl vector is $\rho = 2\La_0 + \varpi$ and the Weyl denominator is
\begin{align*}
A_\rho(q, y) = D(q, y) = y^{1/2} \prod_{n=1}^\infty (1-y^{-1}q^{n-1})(1-q^n)(1-yq^{n}).
\end{align*}
For $\la = \la_0$ and $\la=\la_1$ the numerator of \eqref{eq:adm.char.fla} becomes, respectively, $D(q^5, y q^{-1})$ and $D(q^5, y q^{-2})$. Thus one obtains product formulas
%\begin{align*}
%\chi_{L(\la_0)}(q, y) = {} & \prod_{n=1}^\infty \frac{(1 - y^{-1} q^{5n-4}) (1 - q^{5n}) (1 - y q^{5n-1})}{(1-y^{-1}q^{n-1})(1-q^n)(1-yq^{n})}
%%
%%\chi_{L(\la_1)}(q, y) = {} & \prod_{n=1}^\infty \frac{(1 - y^{-1} q^{5n-3}) (1 - q^{5n}) (1 - y q^{5n-2})}{(1-y^{-1}q^{n-1})(1-q^n)(1-yq^{n})}
%\end{align*}
%and a similar one for $L(\la_1)$.
\begin{align}\label{eq:2/5.sl2.prods}
\chi_{L(\la_0)}(q, y) = \frac{D(q^5, yq^{-1})}{D(q, y)} \quad \text{and} \quad \chi_{L(\la_1)}(q, y) = \frac{D(q^5, yq^{-2})}{D(q, y)}.
\end{align}

\subsection{Affine \texorpdfstring{$W$}{W}-algebras}{\ }

The universal affine vertex algebra of level $k$ associated with finite simple $\g$ is denoted $V^k(\g)$ and its simple quotient $L_k(\g)$. As a $\what\g$-module $V^k(\g)$ is the vacuum module $U(\what\g) \otimes_{U(\what\g_+)} \C \vac$, where $\what\g_+ = \g[t] + \C K + \C d$, acting on $\vac$ as $\g[t] + \C d \mapsto 0$ and $K \mapsto k I$. The vertex algebra structure is generated by quantum fields $a(z) = \sum_{n \in \Z} (at^n) z^{-n-1}$, one for each $a \in \g$. See \cite{FZ92}, \cite[Chapter 4]{va.beginners}.

If $\mbo \subset \g$ is a nilpotent orbit and $f \in \mbo$, the universal affine $W$-algebra $W^k(\g, f)$ is constructed via quantised Drinfeld-Sokolov reduction of $V^k(\g)$ (defined in \cite{FF} for the principal orbit $\mbo$ and \cite{KRW} in general). For $k \neq -h^\vee$ the algebra $W^k(\g, f)$ possesses a conformal structure of central charge \cite[Theorem 2.2]{KRW}
\[
c = \dim(\g_0) - \frac{12}{k+h^\vee} \left|\ov\rho - (k+h^\vee) \frac{h}{2}\right|^2.
\]
The following is a special case of the important structural result \cite[Theorem 4.1]{KW04}.
\begin{thm}\label{thm:PBW}
Let $\{x^1, \ldots x^r\}$ be a basis of $\g^f$, homogeneous with respect to the Dynkin grading. Then $W^k(\g, f)$ contains a strong generating set $\{W^1, \ldots, W^r\}$ where, if $x^i \in \g^f_j$, the conformal weight of $W^i$ is $1-j$. Furthermore PBW monomials in $\{W^1, \ldots, W^r\}$ and their derivatives constitute a basis of $W^k(\g, f)$.
\end{thm}
Since the vertex algebra $W^k(\g, f)$, up to isomorphism, depends on $f$ only through $\mbo$, we often denote it $W^k(\mbo)$, using the Bala-Carter label for $\mbo$. The character of $W^k(\mbo)$ can be determined immediately from Theorem \ref{thm:PBW}, and is given by
\begin{align}\label{eq:univ.W.char}
\chi_{W^k(\mbo)}(q) = \prod_{\substack{j=0 \\ n=1 }}^\infty \frac{1}{(1-q^{j+n})^{\dim(\g^e_j)}}.
\end{align}

The structure of the simple quotient of $W^k(\mbo)$, which we denote $W_k(\mbo)$, is complicated in general. We shall be particularly concerned with cases in which the level $k=-h^\vee + p/u$ is principal admissible. We reserve the notation $W(\mbo, p/u)$ for the simple quotient in this case.

The simple quotient can be obtained by Drinfeld-Sokolov reduction of a specific irreducible highest weight $\what{\g}$-module, more precisely
\[
W(\mbo, p/u) \cong H_{f, -}^0(L(\la)), \quad \la = t_{-h/2}(\phi(\rho)) - \rho.
\]
Here $H_{f, -}^\bullet(-)$ is the ``$-$''-variant of the reduction functor introduced in \cite{FKW}, and $h$ is the neutral element of the $\sll_2$-triple containing $f$. See \cite{FKW} and \cite{Arakawa.rep.I} for $f$ principal and \cite{Arakawa.rep.II} and \cite{AE.JEMS} for other even nilpotent orbits.

The functor $H_{f, -}^\bullet(-)$ is a certain cohomology of the complex $(-) \otimes F(\g_{>0})$. Here $F(\g_{>0})$ is the free fermion vertex superalgebra, it has a strong generating set identified with $\{e_i\} \cup \{\varphi_i\}$, where $\{e_i\}$ is a homogeneous basis of $\g_{>0}$ and $\{\varphi_i\}$ is the dual basis of $\g_{>0}^*$. The character of $W(\mbo, p/u)$ is obtained, via the Euler-Poincar\'{e} principle, as the limit
\begin{align}\label{eq:chi.W.simple}
\chi_{W(\mbo, p/u)} = \lim_{y_i \rightarrow 1} \chi_{L(\la)} \chi_{F(\g_{>0})},
\end{align}
as described in detail in {\cite[Section 2.2]{KW08.Rationality}}. The conformal structure is such that $\chi_{F(\g_{>0})}$ is given by the following variant
\begin{align}\label{eq:denom.qy.prime}
D'(q, y_i) = \prod_{n=1}^\infty \prod_{\al \in \ov\D_{>0}} (1-y_\alpha^{-1} q^{n-1}) (1-y_\alpha q^{n})
\end{align}
of the denominator product \eqref{eq:denom.qy}.

\section{Character formulas}

\subsection{Boundary admissible levels and product formulas}{\ }\label{subsec:bdry}

In general an admissible level for a simple Lie algebra $\g$ of the form $-h^\vee + h^\vee / u$ is said to be boundary admissible. In \cite{KW2017} it was explained how the Kac-Wakimoto character formula naturally leads to product formulas for the characters of boundary admissible $\what\g$-modules, as illustrated in Section \ref{subsec:KW-char-fla} above for the special case $\g = \sll_2$. See also \cite{XYY}.

%If we consider the admissible level $\la = -h^\vee + p/u$ for $p = h^\vee$, then the only possible value of $\nu$ in \eqref{eq:adm.param} is $\nu = \rho$, and the substitution of $D$ for $A_\nu$ in the character formula \eqref{eq:adm.char.fla} yields a product formula for the character of the corresponding irreducible $\what\g$-modules. Such levels, weights and modules, are called boundary admissible.

Combining the product formula \eqref{eq:2/5.sl2.prods} for the character of the $\what\sll_2$-module $L(\la_0)$ of level $-2+2/5$, with equations \eqref{eq:chi.W.simple} and \eqref{eq:denom.qy.prime} for the Drinfeld-Sokolov reduction relative to the standard $\sll_2$-triple, yields
\begin{align}\label{eq:vir-2-5-chars}
\chi_{H^0_{f, -}(L(\la_0))}(q) = {} & \prod_{n=1}^\infty \frac{1}{(1-q^{5n-2})(1-q^{5n-3})}
%
%\chi_{H^0_{f, -}(L(\la_1))}(q) = {} & \prod_{n=1}^\infty \frac{1}{(1-q^{5n-1})(1-q^{5n-4})}
\end{align}
It is well known that the corresponding simple $W$-algebra $W(A_1, 2/5)$ is isomorphic to the Virasoro minimal model $\vir_{2, 5}$, whose character is indeed given by \eqref{eq:vir-2-5-chars}.

%{\color{red}Maybe invent a better notation for $x_0$, or rather for $\beta$, as it is naturally thought of as operating on roots, not coroots.}

The calculation of $\chi_W(q)$ for $W = W(\mbo, h^\vee/u)$ follows the pattern sketched above, and results in the general formula
\begin{align}\label{eq:chi.boundary.fla}
\chi_{W(\mbo, h^\vee/u))}(q) = \lim_{y_i \rightarrow 1} \frac{D(q^u, y_i q^{-\heit(\ov\alpha_i)})}{D(q, y_i)} D'(q, y_i),
\end{align}
where $\heit(\ov\al_i)$ is $j$ for $\ov\al_i \in \g_j$. More explicitly, we introduce the product
\begin{align*}
P_{\mbo, u}(q) = \frac{1}{\prod_{n=1}^\infty (1-q^n)^{\dim(\g_0)}}
\prod_{n=1}^\infty \prod_{j \in \Z} (1-q^{un - j})^{\dim(\g_j)},
\end{align*}
so that after simplification $\chi_{W(\mbo, h^\vee/u))}(q) = P_{\mbo, u}(q)$.
% (Note: really we have also assumed the W-algebra to be lisse, in which case $f \in \mathbb{O}_q$, so $\g_j = 0$ when $|j| \geq u$, guaranteeing that the product does not vanish.)

%To compactly express the infinite products we will be working with, we introduce the notation
%\begin{align*}
%P(q; u; (a_0, \ldots, a_{u-1})) = \prod_{j=0}^{u-1} \prod_{n=1}^\infty (1-q^{un-j})^{-a_j}.
%\end{align*}
%Thus for example we have the following expression for the character of $W = \CW(E_7(a_1), 18/17)$:
%\begin{align*}
%\chi_W(q) = P(q; 17; (0, 0, 1, 1, 1, 1, 2, 2, 2, 2, 2, 2, 1, 1, 1, 1, 0)).
%\end{align*}
%Remarkably we see that the character of $\CW(G_2, 4/17)$ is given by the same infinite product.

For example for $W(G_2, 4/u)$ we read the dimensions $\dim(\g_j)$ off from \eqref{eq:gr.poly.prin} and simplify to obtain its character
\begin{align}\label{eq:G2.17.exact}
P_{G_2, u}(q) = \frac{1}{\prod_{n=1}^\infty (1-q^n)^2} \prod_{\substack{m \in \Z_{> 0} \\ m \equiv \pm 2, \pm 3, \pm 4, \pm 5 \bmod u }} (1-q^{m}).
\end{align}
The product $P_{E_7(a_1), u}(q)$, which is the character of $W(E_7(a_1), 18/u)$, can similarly be read off from \eqref{eq:gr.poly.subreg}. For the particular value $u = 17$ we notice that $P_{E_7(a_1), u}(q)$ and $P_{G_2, u}(q)$ coincide. Therefore isomorphism number (2) of Theorem \ref{thm:isoms.intro} is verified at the level of characters. For reference the series expansion of the product is
\begin{align}\label{eq:G2.E7.boundary.char}
\begin{split}
{} & 1+q^{2}+q^{3}+2 q^{4}+2 q^{5}+5 q^{6}+5 q^{7}+9 q^{8}+11 q^{9}+17 q^{10}+21 q^{11}\\ &+32 q^{12}+39 q^{13}+57 q^{14}+72 q^{15}+99 q^{16}+125 q^{17}+171 q^{18}+214 q^{19} \\ &+286 q^{20}+360 q^{21}+470 q^{22}+590 q^{23}+764 q^{24}+952 q^{25} + \cdots.
\end{split}
\end{align}

\subsection{Near-boundary admissible levels and the Frenkel-Kac isomorphism}{\ }

Suppose $\g$ is simply laced. We recall the Frenkel-Kac isomorphism of vertex algebras
\begin{align*}
L_1(\g) \cong V_Q
\end{align*}
where $Q$ is the root lattice of $\g$ and $V_Q$ is the corresponding lattice vertex algebra (see \cite{Frenkel-Kac} and \cite[Chapter 5]{va.beginners}). From the construction of $V_Q$ one has the character formula
\begin{align}\label{eq:lattice.char}
\chi_{L(\La_0)}(q, y) = \chi_{V_Q}(q, y) = \sum_{\alpha \in Q} y^{\alpha} q^{|\alpha|^2/2},
\end{align}
where $y$ is shorthand for $(y_1, \ldots, y_\ell)$ and $y^\alpha$ is shorthand for $y_i^{m_1} \cdots y_\ell^{m_\ell}$ where $m_1, \ldots m_\ell \in \Z$ are the coefficients of $\al$ relative to the basis $\{\al_1, \ldots, \al_\ell\}$. For definiteness we write $\Theta_Q(q, y)$ for the sum on the right hand side of \eqref{eq:lattice.char}. Then the Weyl denominator formula \eqref{eq:denom.fla} yields
\begin{align}\label{eq:A.at.hvplus1}
A_{(h^\vee+1)\La_0 + \ov\rho}(q, y) = D(q, y) \Theta_Q(q, y).
\end{align}
Equations \eqref{eq:lattice.char} and \eqref{eq:A.at.hvplus1} provide an alternative to the computation of $A_{(h^\vee+1)\La_0 + \ov\rho}$ from its definition as a sum over the Weyl group. This is particularly convenient if the Weyl group of $\g$ is very large. We may apply equation \eqref{eq:A.at.hvplus1} to obtain character formulas for simple $W$-algebra at ``near boundary'' admissible level, by which we mean levels of the form $k = -h^\vee + p/u$ for $p = h^\vee+1$. The character of $W(\mbo, (h^\vee+1)/u)$ is obtained in the same way as in Subsection \ref{subsec:bdry}, except that in formula \eqref{eq:chi.boundary.fla} the numerator $D$ is replaced with $D \cdot \Theta_Q$. The change of variable $y_i \mapsto y_i q^{-\heit(\ov\al_i)}$ in the argument of $\Theta_Q$ transforms it to
\begin{align*}
\sum_{\al \in Q} q^{u |\al|^2/2 - (\xi, \al)},
\end{align*}
where $\xi = \sum_{i=1}^\ell \heit(\ov\al_i) \varpi_i$.

%%%%%%%%%%%%%%%%%%%%%%%%%%%%%%%%
\subsubsection{Characters of $W(E_8(a_1), 31/u)$}{\ }\label{subsubsec:E8.near.admissible}

If $\g = E_8$ then the root lattice $Q(E_8)$ may be presented as the subset of
\begin{align*}
\Z^8 \cup (1/2+\Z)^8 \subset \mathbb{R}^8
\end{align*}
consisting of vectors, the sum of whose entries is an even integer.
%A trivial observation that nevertheless helped me in the Mathematica computations is the following: the set $((1/2 + \Z)^8)^{\text{even}}$ coincides with the set $\{(1/2+m_1, \ldots, 1/2+m_8) | \text{$\sum_i m_i$ is even} \}$.
A choice of simple roots is given by the rows of the matrix
\begin{align*}
\left[\begin{array}{cccccccc}
 1 & -1 & 0 & 0 & 0 & 0 & 0 & 0 \\
 0 &  1 & -1 & 0 & 0 & 0 & 0 & 0 \\
 0 & 0 &  1 & -1 & 0 & 0 & 0 & 0 \\
 0 & 0 & 0 &  1 & -1 & 0 & 0 & 0 \\
 0 & 0 & 0 & 0 &  1 & -1 & 0 & 0 \\
 0 & 0 & 0 & 0 & 0 &  1 & -1 & 0 \\
 0 & 0 & 0 & 0 & 0 & 1 &  1 & 0 \\
 -1/2 & -1/2 & -1/2 & -1/2 & -1/2 & -1/2 & -1/2 & -1/2 \\
      \end{array}
\right].
\end{align*}

If we restrict attention to the case of subregular nilpotent $\mbo = E_8(a_1)$, then $\heit(\al_{*}) = 0$ for $*$ the trivalent node of the Dynkin diagram of $E_8$ and $\heit(\al_i) = 2$ for the other simple roots. The corresponding vector is $\xi_8 = (5, 4, 3, 2, 1, 1, 0, -18)$. It is then straightforward to compute the character $\chi_W(q)$ of $W = W(E_8(a_1), 31/u)$ to high order in $q$. For $25 \leq u \leq 29$ we obtain
% \begin{align*}
% &u = 25: & \chi_W(q) &= 1+q^{2}+q^{3}+q^{4}+q^{5}+2 q^{6}+2 q^{7}+3 q^{8}+3 q^{9}+4 q^{10}+4 q^{11}+6 q^{12}+6 q^{13}+8 q^{14}+9 q^{15}+11 q^{16}+12 q^{17}+15 q^{18}+16 q^{19}+20 q^{20}+22 q^{21}+26 q^{22}+29 q^{23}+35 q^{24}+38 q^{25} + \cdots \\
% %
% &u = 26: & \chi_W(q) &= 1+q^{2}+q^{3}+2 q^{4}+2 q^{5}+4 q^{6}+4 q^{7}+7 q^{8}+8 q^{9}+12 q^{10}+14 q^{11}+20 q^{12}+23 q^{13}+32 q^{14}+39 q^{15}+51 q^{16}+61 q^{17}+80 q^{18}+95 q^{19}+122 q^{20}+146 q^{21}+183 q^{22}+219 q^{23}+273 q^{24}+324 q^{25} + \cdots \\
% %
% &u = 27: & \chi_W(q) &= 1+q^{2}+q^{3}+2 q^{4}+2 q^{5}+5 q^{6}+5 q^{7}+9 q^{8}+11 q^{9}+17 q^{10}+21 q^{11}+32 q^{12}+39 q^{13}+56 q^{14}+71 q^{15}+97 q^{16}+122 q^{17}+166 q^{18}+207 q^{19}+275 q^{20}+345 q^{21}+448 q^{22}+560 q^{23}+722 q^{24}+896 q^{25} + \cdots \\
% %
% &u = 28: & \chi_W(q) &= 1+q^{2}+q^{3}+2 q^{4}+2 q^{5}+5 q^{6}+5 q^{7}+10 q^{8}+12 q^{9}+19 q^{10}+24 q^{11}+38 q^{12}+47 q^{13}+70 q^{14}+91 q^{15}+128 q^{16}+165 q^{17}+230 q^{18}+295 q^{19}+402 q^{20}+518 q^{21}+690 q^{22}+886 q^{23}+1170 q^{24}+1491 q^{25} + \cdots \\
% %
% &u = 29: & \chi_W(q) &= 1+q^{2}+q^{3}+2 q^{4}+2 q^{5}+5 q^{6}+5 q^{7}+10 q^{8}+12 q^{9}+20 q^{10}+25 q^{11}+40 q^{12}+50 q^{13}+76 q^{14}+100 q^{15}+143 q^{16}+187 q^{17}+265 q^{18}+345 q^{19}+478 q^{20}+625 q^{21}+847 q^{22}+1105 q^{23}+1483 q^{24}+1922 q^{25} + \cdots.
% \end{align*}
\begin{align*}
&u = 25: & \chi_W(q) = {} & 1+q^{2}+q^{3}+q^{4}+q^{5}+2 q^{6}+2 q^{7}+3 q^{8}+3 q^{9}+4 q^{10}+4 q^{11}\\&& &+6 q^{12}+6 q^{13}+8 q^{14}+9 q^{15}+11 q^{16}+12 q^{17}+15 q^{18}+16 q^{19}\\&& &+20 q^{20}+22 q^{21}+26 q^{22}+29 q^{23}+35 q^{24}+38 q^{25} + \cdots \\
&u = 26: & \chi_W(q) = {} & 1+q^{2}+q^{3}+2 q^{4}+2 q^{5}+4 q^{6}+4 q^{7}+7 q^{8}+8 q^{9}+12 q^{10}+14 q^{11}\\&& &+20 q^{12}+23 q^{13}+32 q^{14}+39 q^{15}+51 q^{16}+61 q^{17}+80 q^{18}+95 q^{19}\\&& &+122 q^{20}+146 q^{21}+183 q^{22}+219 q^{23}+273 q^{24}+324 q^{25} + \cdots \\
&u = 27: & \chi_W(q) = {} & 1+q^{2}+q^{3}+2 q^{4}+2 q^{5}+5 q^{6}+5 q^{7}+9 q^{8}+11 q^{9}+17 q^{10}+21 q^{11}\\&& &+32 q^{12}+39 q^{13}+56 q^{14}+71 q^{15}+97 q^{16}+122 q^{17}+166 q^{18}+207 q^{19}\\&& &+275 q^{20}+345 q^{21}+448 q^{22}+560 q^{23}+722 q^{24}+896 q^{25} + \cdots \\
&u = 28: & \chi_W(q) = {} & 1+q^{2}+q^{3}+2 q^{4}+2 q^{5}+5 q^{6}+5 q^{7}+10 q^{8}+12 q^{9}+19 q^{10}+24 q^{11}\\&& &+38 q^{12}+47 q^{13}+70 q^{14}+91 q^{15}+128 q^{16}+165 q^{17}+230 q^{18}+295 q^{19}\\&& &+402 q^{20}+518 q^{21}+690 q^{22}+886 q^{23}+1170 q^{24}+1491 q^{25} + \cdots \\
&u = 29: & \chi_W(q) = {} & 1+q^{2}+q^{3}+2 q^{4}+2 q^{5}+5 q^{6}+5 q^{7}+10 q^{8}+12 q^{9}+20 q^{10}+25 q^{11}\\&& &+40 q^{12}+50 q^{13}+76 q^{14}+100 q^{15}+143 q^{16}+187 q^{17}+265 q^{18}+345 q^{19}\\&& &+478 q^{20}+625 q^{21}+847 q^{22}+1105 q^{23}+1483 q^{24}+1922 q^{25} + \cdots.
\end{align*}
\begin{rem}
In \cite{AE.JEMS} it was proved that
\begin{align*}
W(E_8(a_1), 31/25) \cong \vir_{2, 5} \quad \text{and} \quad W(E_8(a_1), 31/26) \cong \vir_{2, 13}.
\end{align*}
From the first of these isomorphisms we deduce the following pheasant identity
\begin{align*}
\sum_{\al \in Q(E_8)} q^{25 |\al|^2/2 - (\xi_8, \al)} = \prod_{n=1}^\infty \frac{(1-q^{25n-10})(1-q^{25n-15})}{(1-q^{5n-2})(1-q^{5n-3})}.
\end{align*}
\end{rem}

\subsubsection{Character of $W(E_7(a_1), 19/16)$}{\ }

The root lattice $Q(E_7)$ may be realised as the subset of $Q(E_8)$ consisting of vectors, the sum of whose entries vanishes, i.e., as the intersection of $Q(E_8)$ with the orthogonal complement of $(1, 1, 1, 1, 1, 1, 1, 1)$. As for $E_8$, the subregular Dynkin grading assigns $\heit(\alpha_{*}) = 0$ to the trivalent node $*$, and $\heit(\al_i) = 2$ to the other simple roots. A convenient choice of simple roots in $Q(E_7)$ is taken, and the corresponding vector $\xi$ is computed. A short calculation then yields $\xi = \xi_7 = (8, -5, -4, -3, -2, -1, 0)$. To prove isomorphism number (1) of Theorem \ref{thm:isoms.intro}, we set $u=16$ to obtain the character of $W = W(E_7(a_1), 19/16)$ as
\begin{align}\label{eq:E7a1.19.16.char}
\begin{split}
\chi_W(q) {} &= 1+q^{2}+q^{3}+3 q^{4}+3 q^{5}+7 q^{6}+8 q^{7}+14 q^{8}+18 q^{9}+28 q^{10}+36 q^{11} \\
&+55 q^{12} +70 q^{13}+101 q^{14}+131 q^{15}+182 q^{16}+234 q^{17}+319 q^{18}+408 q^{19} \\
&+544 q^{20} +694 q^{21} +909 q^{22}+1153 q^{23}+1494 q^{24}+1881 q^{25} + \cdots.
\end{split}
\end{align}

\subsection{Other character computations}{\ }

In this subsection we compute series expansions of characters of the principal $W$-algebras of types $G_2$, $B_3$ and $F_4$ appearing in the right hand side of \eqref{eq:isoms.intro}. Case (2) is a boundary level for $G_2$ and so we have the exact expression for the character given in \eqref{eq:G2.17.exact} above. For the others we compute directly by summing over the Weyl group. In organising the computations we follow the conventions of Bourbaki {\cite[pp. 265-290]{Bourbaki.Lie.4-6}}: for $G_2$ the short simple root is denoted $\al_1$ and the long one $\al_2$, etc.

\subsubsection{Character of $W(G_2, 9/7)$}{\ }

As in Subsection \ref{subsec:KW-char-fla} above, the character of $W = W(G_2, p/u)$ is given by
\begin{align}\label{eq:W.GBF.prin.general}
\chi_{W}(q) = \frac{A_{p\La_0 + \ov\rho}(q^u, q^{-1}, q^{-1})}{\prod_{n=1}^\infty (1-q^n)^\ell},
\end{align}
where $\ell = 2$ is the rank and $\ov\rho$ is the finite Weyl vector of $G_2$. The sum over the affine Weyl group defining $A_{p\La_0 + \ov\rho}$ may be expressed as a sum over the finite Weyl group of certain theta functions associated with the coroot lattice $Q^\vee$, namely
\begin{align}\label{eq:A-from-T.GBF}
A_{p\La_0 + \ov\rho}(q, y_1, \ldots y_\ell) = \sum_{w \in W} (-1)^w \Theta^{(p)}_{w(\ov\rho)}(q, y_1, \ldots, y_\ell)
\end{align}
where
\begin{align}\label{eq:T-from-S.GBF}
\Theta^{(p)}_{\mu}(q, y_1, \ldots, y_\ell) = \sum_{\al \in Q^\vee} q^{(\mu, \alpha) + p (\al, \al)/2} e^{2\pi i (z, \mu + p \alpha)}.
\end{align}
More explicitly, we have $Q^\vee = \Z \cdot 3\al_1 + \Z \cdot \al_2$, and if we write $\al = 3 m_1 \al_1 + m_2 \al_2$ and $\mu = c_1 \al_1 + c_2 \al_2$ then $e^{2\pi i (z, \mu + p \alpha)}$ becomes $y_1^{c_1+3p m_1} y_2^{c_2+p m_2}$ and the exponent of $q$ becomes $(c_1, c_2)M \binom{3m_1}{m_2} + \tfrac{1}{2} p (3m_1, m_2)M \binom{3m_1}{m_2}$ where $M$ is the Gram matrix
\begin{align*}
M = \left[\begin{array}{cc}
2/3 & -1 \\
-1 & 2 \\
\end{array}\right].
\end{align*}
Relative to the basis $\{\al_1, \al_2\}$ the finite Weyl group is generated by the reflections
\begin{align*}
r_1 = \left[\begin{array}{cc}
-1 & 3 \\
0 & 1 \\
\end{array}\right] \quad \text{and} \quad
r_2 = \left[\begin{array}{cc}
1 & 0 \\
1 & -1 \\
\end{array}\right].
\end{align*}
The images of $\ov\rho = 5\alpha_1 + 3\alpha_2$ under the finite Weyl group are computed easily enough to be
\begin{align*}
&{}^+(5, 3) &  &{}^-(4, 3) &  &{}^-(5, 2) &  &{}^+(1, 2) &  &{}^+(4, 1) &  &{}^-(-1, 1) \\
&{}^-(1, -1) &  &{}^+(-4, -1) &  &{}^+(-1, -2) &  &{}^-(-5, -2) &  &{}^-(-4, -3) &  &{}^+(-5, -3),
\end{align*}
(the superscripts recording $(-1)^w$). Substitution into \eqref{eq:W.GBF.prin.general} with $p/u = 9/7$ yields the series expansion
\begin{align*}%\label{eq:W.G2.prin.9/7}
\chi_{W}(q) = {} & 1 + q^2 + q^3 + 2 q^4 + 2 q^5 + 5 q^6 + 5 q^7 + 9 q^8 + 11 q^9 +  17 q^{10} + 21 q^{11} \\
&+ 32 q^{12} + 39 q^{13} + 56 q^{14} + 71 q^{15} +  97 q^{16} + 122 q^{17} + 166 q^{18} + 207 q^{19} \\
&+ 275 q^{20} + 345 q^{21} +  448 q^{22} + 560 q^{23} + 722 q^{24} + 896 q^{25} + \cdots
\end{align*}
%I'm very pleased to say that this agrees with the character of subregular near-boundary $E_8$ at denominator $27$.

\subsubsection{Character of $W(B_3, 8/7)$}{\ }

%As in Subsection \ref{subsec:KW-char-fla} above, the character of $W = W(B_3, p/u)$ is given by
%\begin{align}\label{eq:W.B3.prin.general}
%\chi_{W}(q) = \frac{A_{p\La_0 + \ov\rho}(q^u, q^{-1}, q^{-1})}{\prod_{n=1}^\infty (1-q^n)^3}.
%\end{align}
%The sum over the affine Weyl group defining $A_{p\La_0 + \ov\rho}(q, y_1, y_2, y_3)$ may be expressed as a sum over the finite Weyl group of certain theta functions associated with the coroot lattice $Q^\vee$. In this case $Q^\vee = \Z \cdot \al_1 + \Z \cdot \al_2 + \Z 2 \al_3$. Unwinding the definition of $A$ yields
%\begin{align*}
%A_{p\La_0 + \ov\rho}(q, y_1, y_2, y_3) = \sum_{w \in W} (-1)^w \Theta^{(p)}_{w(\ov\rho)}(q, y_1, y_2, y_3)
%\end{align*}
%where
%\begin{align*}
%\Theta^{(p)}_{\mu}(q, y_1, y_2, y_3) = \sum_{\al \in Q^\vee} q^{(\mu, \alpha) + p (\al, \al)/2} e^{2\pi i (z, \mu + p \alpha)}.
%\end{align*}

We have \eqref{eq:A-from-T.GBF} and \eqref{eq:T-from-S.GBF} as before, and the character $\chi_W(q)$ of $W = W(B_3, p/u)$ is given by \eqref{eq:W.GBF.prin.general}, having substituted $y_1 = y_2 = y_3 = q^{-1}$.

More explicitly, we have $Q^\vee = \Z \cdot \al_1 + \Z \cdot \al_2 + \Z \cdot 2 \al_3$, and if we write $\al = m_1 \al_1 + m_2 \al_2 + 2m_3 \al_3$ and $\mu = c_1 \al_1 + c_2 \al_2 + c_3\al_3$ then $e^{2\pi i (z, \mu + p \alpha)}$ becomes $y_1^{c_1+p m_1} y_2^{c_2+p m_2} y_3^{c_3+2p m_3}$ and the exponent of $q$ becomes $(c_1, c_2, c_3) M (m_1, m_2, 2m_3)^T + \tfrac{1}{2} p (m_1, m_2, 2m_3) M  (m_1, m_2, 2m_3)^T$ where $M$ is the Gram matrix
\begin{align*}
M = \left[\begin{array}{ccc}
2 & -1 & 0 \\
-1 & 2 & -1 \\
0 & -1 & 1 \\
\end{array}\right].
\end{align*}
The finite Weyl group has cardinality $48$ and, relative to the basis $\{\al_1, \al_2, \al_3\}$, is generated by the reflection matrices
\begin{align*}
r_1 = \left[\begin{array}{ccc}
-1 & 1 & 0 \\
0 & 1 & 0 \\
0 & 0 & 1 \\
\end{array}\right], \quad
r_2 = \left[\begin{array}{ccc}
1 & 0 & 0 \\
1 & -1 & 1 \\
0 & 0 & 1 \\
\end{array}\right] \quad \text{and} \quad
r_3 = \left[\begin{array}{ccc}
1 & 0 & 0 \\
0 & 1 & 0 \\
0 & 2 & -1 \\
\end{array}\right].
\end{align*}
The finite Weyl vector is $\ov\rho = \frac{1}{2}\left( 5\alpha_1 + 8 \alpha_2 + 9 \alpha_3 \right)$ and, upon computing its orbit under the finite Weyl group and substituting into \eqref{eq:W.GBF.prin.general} with $p/u = 8/7$, we obtain the series expansion
\begin{align}\label{eq:B3.simple.char}
\begin{split}
\chi_{W}(q) = {} & 1+q^2+q^3+3 q^4+3 q^5+7 q^6+8 q^7+14 q^8+18 q^9+28 q^{10}+36 q^{11} \\
&+55 q^{12} +70 q^{13}+101 q^{14}+131 q^{15}+182 q^{16}+234 q^{17}+319 q^{18}+408 q^{19} \\
&+544 q^{20}+694 q^{21} +909 q^{22} +1153 q^{23}+1494 q^{24}+1881 q^{25} + \cdots
\end{split}
\end{align}

\subsubsection{Character of $W(F_4, 14/13)$}{\ }

We have \eqref{eq:A-from-T.GBF} and \eqref{eq:T-from-S.GBF} as before, and the character $\chi_W(q)$ of $W = W(F_4, p/u)$ is given by \eqref{eq:W.GBF.prin.general}, having substituted $y_1 = y_2 = y_3 = y_4 = q^{-1}$.

More explicitly, we have $Q^\vee = \Z \cdot \al_1 + \Z \cdot \al_2 + \Z \cdot 2\al_3 + \Z \cdot 2\al_4$, and if we write $\al = m_1 \al_1 + m_2 \al_2 + 2m_3 \al_3 + 2m_4 \al_4$ and $\mu = c_1 \al_1 + c_2 \al_2 + c_3\al_3 + c_4 \al_4$ then $e^{2\pi i (z, \mu + p \alpha)}$ becomes $y_1^{c_1+p m_1} y_2^{c_2+p m_2} y_3^{c_3+2p m_3} y_4^{c_4+2p m_4}$ and the exponent of $q$ becomes
\[
(c_1, c_2, c_3, c_4) M (m_1, m_2, 2m_3, 2m_4)^T + \tfrac{1}{2} p (m_1, m_2, 2m_3, 2m_4) M  (m_1, m_2, 2m_3, 2m_4)^T
\]
where $M$ is the Gram matrix
\begin{align*}
M = \left[\begin{array}{cccc}
2 & -1 & 0 & 0 \\
-1 & 2 & -1 & 0 \\
0 & -1 & 1 & -1/2 \\
0 & 0 & -1/2 & 1 \\
\end{array}\right].
\end{align*}
The finite Weyl group has cardinality $1152$ and, relative to the basis $\{\al_1, \al_2, \al_3, \al_4\}$, is generated by the reflection matrices
{\tiny
\begin{align*}
r_1 = \left[\begin{array}{cccc}
-1 & 1 & 0 & 0 \\
0 & 1 & 0 & 0 \\
0 & 0 & 1 & 0 \\
0 & 0 & 0 & 1 \\
\end{array}\right], \quad
r_2 = \left[\begin{array}{cccc}
1 & 0 & 0 & 0 \\
1 & -1 & 1 & 0 \\
0 & 0 & 1 & 0 \\
0 & 0 & 0 & 1 \\
\end{array}\right], \quad
r_3 = \left[\begin{array}{cccc}
1 & 0 & 0 & 0 \\
0 & 1 & 0 & 0 \\
0 & 2 & -1 & 1 \\
0 & 0 & 0 & 1 \\
\end{array}\right], \quad \text{and} \quad
r_4 = \left[\begin{array}{cccc}
1 & 0 & 0 & 0 \\
0 & 1 & 0 & 0 \\
0 & 0 & 1 & 0 \\
0 & 0 & 1 & -1 \\
\end{array}\right].
\end{align*}
}
The finite Weyl vector is $\ov\rho = 8 \alpha_1 + 15 \alpha_2 + 21 \alpha_3 + 11 \alpha_4$ and, upon computing its orbit under the finite Weyl group and substituting into \eqref{eq:W.GBF.prin.general} with $p/u = 14/13$, we obtain the series expansion
\begin{align*}%\label{eq:W.G2.prin.9/7}
\chi_{W}(q) = {} & 1 + q^2 + q^3 + 2 q^4 + 2 q^5 + 5 q^6 + 5 q^7 + 10 q^8 + 12 q^9 + 19 q^{10} + 24 q^{11} \\
&+ 38 q^{12} + 47 q^{13} + 70 q^{14} + 91 q^{15} + 128 q^{16} + 165 q^{17} + 230 q^{18} + 295 q^{19} \\
&+ 402 q^{20} + 518 q^{21} + 690 q^{22} + 886 q^{23} + 1170 q^{24} + 1491 q^{25} + \cdots
\end{align*}

\section{Isomorphisms}\label{sec:isoms}

In this section we prove the isomorphisms \eqref{eq:isoms.intro} using the character computations of the preceding sections, explicit computations with OPE, and well-known results of the Virasoro Lie algebra from Section \ref{sec: Representation theory of the Virasoro Lie algebra}.

%%%%%%%%%%%%%%%%%%%%%%%%%%%
%%%%%%%%%%%%%%%%%%%%%%%%%%%
% CASE 2
%%%%%%%%%%%%%%%%%%%%%%%%%%%
%%%%%%%%%%%%%%%%%%%%%%%%%%%

\subsection{Proof of isomorphism number (2)}{\ }

We write $k = -h^\vee_{G_2} + 4/17$, so that $W^k(G_2)$ has central charge $c = -1420/17$. We shall prove the existence of a surjective homomorphism
\[
W^k(G_2) \rightarrow W(E_7(a_1), 18/17).
\]
The isomorphism asserted in Theorem \ref{thm:isoms.intro} follows immediately from simplicity of the target.

From equation \eqref{eq:univ.W.char} we have
\begin{align}\label{eq:univ.G2.char}
\begin{split}
\chi_{W^k(G_2)}(q) = {} & 1 + q^2 + q^3 + 2 q^4 + 2 q^5 + 5 q^6 + 5 q^7 + 9 q^8 + 11 q^9 + 17 q^{10} \\
&+ 21 q^{11} + 33 q^{12} + 40 q^{13} + \cdots
\end{split}
\end{align}
For the rest of this subsection we write $W$ for $W(E_7(a_1), 18/17)$ and $\wtil{W}$ for the corresponding universal $W$-algebra. By Theorem \ref{thm:PBW} the vertex algebra $\wtil{W}$ has a strong generating set with conformal weights
\begin{align}\label{eq:E7a1.gens}
2, 4, 6, 6, 8, 9, 10, 12, 14,
\end{align}
and has a basis consisting of PBW monomials in these generators and their derivatives. We choose such a set of generators, and denote their images in the simple quotient $W$ by $L$, $W^4$, $W^{6, 1}$, $W^{6, 2}$, $W^8$, $W^9$, $W^{10}$, $W^{12}$ and $W^{14}$, respectively.

We recall $\chi_W(q) = \chi_{W(G_2, 4/17)}(q)$, with both sides given by the series \eqref{eq:G2.E7.boundary.char} which, by inspection, coincides with \eqref{eq:univ.G2.char} above modulo $q^{12}$.

Let us denote by $M^0$ the $U(\CL)$-submodule of $W$ generated by $\vac$, and $M^1$ the $U(\CL)$-submodule generated by $\vac$ and $W^4$. Since $W_{\leq 3} = M^0_{\leq 3}$ we have $L_{n} W^4 \in M^0$ for all $n \geq 1$, and we obtain a short exact sequence of $U(\CL)$-modules
\begin{align}\label{eq:SES}
0 \rightarrow M^0 \rightarrow M^1 \rightarrow Q \rightarrow 0,
\end{align}
in which $M^0$ is a quotient of $M(c, 0)$ and $Q$ is a quotient of $M(c, 4)$. In fact $M^0$ is a nonzero quotient of the vacuum $U(\CL)$-module $V(c)$, and $M^0 \cong V(c)$ since the latter is irreducible. By comparing $q$-characters up to $q^4$, we see that $Q = 0$ and hence $M^1 = M^0$.

We now repeat the argument with $M^0 \subset W$ as before, and $M^1 \subset W$ the $U(\CL)$-submodule generated by $\vac$, $W^{6, 1}$ and $W^{6, 2}$. Once again we have a short exact sequence \eqref{eq:SES}, and by comparing $q$-characters we see that $\dim Q_{6} = 1$ and so $Q$ is a quotient of $M(c, 6)$. In fact $Q \cong M(c, 6)$ since the latter is irreducible. We choose a generator of $Q$ of conformal weight $6$ and denote by $W^6$ a preimage of this generator in $W$. In fact we have $M^1 \cong V(c) \oplus M(c, 6)$ since there are no nontrivial extensions between the components, so without loss of generality we may take $W^6$ to be a primary vector.

We may continue in the same way, using the equality
\begin{align*}
\chi_{W}(q) \equiv \chi_{\wtil{W}^k(G_2)}(q) \quad \bmod{q^{12}}
\end{align*}
to deduce that $W_{\leq 11} = M_{\leq 11}$ where now $M \subset W$ denotes the $U(\CL)$-submodule generated by $\vac$ and $W^6$.

Now we study the possible OPE relations between $L$ and the primary vector $W^6$. In general
\begin{align*}
[L_\la W^6] = TW^6 + 6\la W^6, \qquad
[{W^6}_{\la}W^6] = \sum_{k=0}^{11} A_{k} \lambda^{11-k},
\end{align*}
for some $A_k \in V_k$. We write $A_k$ relative to a PBW monomial basis of $(V(c) \oplus M(c, 6))_k$, expand the commutator formula \eqref{eq:comm.fla} for the triples $(L, W^6, W^6)$ and $(W^6, W^6, W^6)$ relative to the chosen basis, and equate coefficients. Thus we obtain a system of quadratic equations (most of them linear in fact) on the $75$ coefficients describing the $A_k$ for $0 \leq k \leq 11$. We find that this system has a unique $1$-parameter family of solutions with parameter given by the coefficient $a_0\in \C$ in $A_0=a_0 \vac$. The simplicity of $W$ implies $a_0\neq0$.
Since we have not yet normalised $W^6$, we may fix such a normalisation to obtain  a unique solution so that $a_0=1$.

The unique solution must describe both $W^k(G_2)$ and $W$, and it follows from the universal property of the former that there exists a homomorphism of vertex algebras
\begin{align*}
\varphi : W^k(G_2) \rightarrow W.
\end{align*}
From general principles we have $\chi_{W^k(G_2)}(q) \geq \chi_{\im(\varphi)}(q) \geq \chi_{W_k(G_2)}(q)$ and $\chi_{W}(q) \geq \chi_{\im(\varphi)}(q)$, but since we have also computed $\chi_{W}(q) = \chi_{W_k(G_2)}(q)$ we in fact have
\[
W = \im(\varphi) \cong W_k(G_2) = W(G_2, 4/17).
\]

\subsection{Proof of isomorphism number (3)}{\ }

The analysis of this case is similar to the previous one. We shall prove the existence of a surjective homomorphism
\[
W^k(G_2) \rightarrow W = W(E_8(a_1), 31/27),
\]
where now $c = -590/9$. By Theorem \ref{thm:PBW} the universal vertex algebra $\wtil{W} = W^k(E_8(a_1))$ has a set of strong generators with conformal weights
\[
2, 6, 8, 10, 12, 14, 15, 18, 20, 24.
\]
We have the formula in Section \ref{subsubsec:E8.near.admissible}, at $u=27$, for the character. We consider the $U(\CL)$-module structure of $W$ and deduce, just as before, that $W_{\leq 11} = M_{\leq 11}$ where $M$ is the $U(\CL)$-submodule of $W$ generated by $\vac$ and a primary vector $W^6$ of conformal weight $6$.

It follows that there exists a homomorphism of vertex algebras
\begin{align*}
\varphi : W^k(G_2) \rightarrow W,
\end{align*}
and we have $\chi_{W^k(G_2)}(q) \geq \chi_{\im(\varphi)}(q) \geq \chi_{W_k(G_2)}(q)$ and $\chi_{W}(q) \geq \chi_{\im(\varphi)}(q)$. We have not (yet) established $\chi_{W}(q) = \chi_{W_k(G_2)}(q)$ as we had in the previous case, but we have computed explicitly $\chi_{W}(q) \equiv \chi_{W_k(G_2)}(q) \bmod{q^{25}}$. From this we deduce that the generators $W^{14}, \ldots, W^{24}$ must lie in $\im(\varphi)$. Therefore the conclusion
\[
W = \im(\varphi) \cong W_k(G_2) = W(G_2, 9/7)
\]
follows as in the previous case.

%%%%%%%%%%%%%%%%%%%%%%%%%%%
%%%%%%%%%%%%%%%%%%%%%%%%%%%
% CASE 1
%%%%%%%%%%%%%%%%%%%%%%%%%%%
%%%%%%%%%%%%%%%%%%%%%%%%%%%

\subsection{Proof of isomorphism number (1)}{\ }

The central charge of both vertex algebras in this case is $c = -135/8$. We shall prove the existence of a surjective homomorphism
\[
W^k(B_3) \rightarrow W = W(E_7(a_1), 19/16).
\]
From equation \eqref{eq:univ.W.char} we have
\begin{align*}
\chi_{W^k(B_3)}(q) = {} & 1 + q^2 + q^3 + 3 q^4 + 3 q^5 + 7 q^6 + 8 q^7 + 15 q^8 + 19 q^9 + 32 q^{10} \\
&+ 42 q^{11} + 68 q^{12} + 89 q^{13} + \cdots
\end{align*}
We choose a strong generating set of the universal vertex algebra $\wtil{W} = W^k(E_7(a_1))$ with conformal weights as in \eqref{eq:E7a1.gens}. We have the formulas \eqref{eq:E7a1.19.16.char} and \eqref{eq:B3.simple.char} for the characters of $W$ and $W(B_3, 8/7)$, which we see coincide modulo $q^{25}$.

Let us denote by $M^0$ the $U(\CL)$-submodule of $W$ generated by $\vac$, and $M^1$ the $U(\CL)$-submodule generated by $\vac$ and $W^4$. Since $W_{\leq 3} = M^0_{\leq 3}$ we have $L_{n} W^4 \in M^0$ for all $n \geq 1$, and we obtain a short exact sequence of $U(\CL)$-modules
\begin{align}\label{eq:case.1.SES.1}
0 \rightarrow M^0 \rightarrow M^1 \rightarrow Q \rightarrow 0,
\end{align}
in which $M^0$ is a quotient of $M(c, 0)$ and $Q$ is a quotient of $M(c, 4)$. In fact $M^0$ is a nonzero quotient of the vacuum $U(\CL)$-module $V(c)$, and $M^0 \cong V(c)$ since the latter is irreducible. By comparing $q$-characters, we see that $Q \neq 0$ and hence $W^4 \neq 0$. We also note that $Q \cong M(c, 4)$ as the latter is irreducible, and so we know the character of $M^1$.

We now repeat the argument with $M^1 \subset W$ as before, and $M^2 \subset W$ the $U(\CL)$-submodule generated by $\vac$, $W^4$, $W^{6, 1}$ and $W^{6, 2}$. Once again we have a short exact sequence
\begin{align}\label{eq:case.1.SES.2}
0 \rightarrow M^1 \rightarrow M^2 \rightarrow Q \rightarrow 0,
\end{align}
and by comparing $q$-characters we see that $\dim Q_{6} = 1$ and so $Q$ is a quotient of $M(c, 6)$. In fact $Q \cong M(c, 6)$ since the latter is irreducible. We choose a generator of $Q$ of conformal weight $6$ and denote by $W^6$ a preimage of this generator in $W$. In fact we have $M^2 \cong V(c) \oplus M(c, 4) \oplus M(c, 6)$ since there are no nontrivial extensions between the components, and we may take the generators $W^4$ and $W^6$ to be primary. We have the character
\begin{align*}
\chi_{M^2}(q) = {} & 1 + q^2 + q^3 + 3 q^4 + 3 q^5 + 7 q^6 + 8 q^7 + 14 q^8 + 18 q^9 + 28 q^{10} \\
&+ 36 q^{11} + 54 q^{12} + 69 q^{13} + 98 q^{14} + 127 q^{15} + 174 q^{16} + 223 q^{17} + \cdots
\end{align*}

We continue, defining $M^3 \subset W$ to be the $U(\CL)$-submodule generated by $\vac$, $W^4$, $W^{6}$ and $W^8$. Now since we have equality $\chi_{W}(q) \equiv \chi_{M^2}(q) \bmod{q^{9}}$, in fact $M^3 = M^2$ and so $W^8 \in M^2$. Continuing in this way we find that $W^9, W^{10} \in M^2$ also. We also conclude that $W_{\leq 11} = M_{\leq 11}$ where now $M \subset W$ denotes the $U(\CL)$-submodule generated by $\vac$, $W^4$ and $W^6$. In particular all coefficients of $[W^4_{\la}W^4]$, $[W^4_{\la}W^6]$ and $[W^6_{\la}W^6]$ lie within $M$.

%Now we come to conformal weight $12$, where we find $\dim(W_{12}) = 55$ while $\dim(M_{12})=54$. Thus there are two possbilities: (1) $W^{12} = 0$ and (2) $W^{12} \neq 0$ and some singular vector appears in the $<2, 4, 6>$-generated algebra.
%
%{\color{red}We would like to RULE OUT the second case.
%
%NO, WAIT... Something is weird. Solving the OPE equations with all the }

Now we study the possible OPE relations between $L$ and the primary vectors $W^4$ and $W^6$. We posit
% \begin{align}\label{eq:generalOPE.case1}
% \begin{split}
% [{W^4}_{\la}W^4] &= \sum_{k=0}^{7} A_{k} \lambda^{7-k}, \\
% %
% [{W^4}_{\la}W^6] &= \sum_{k=0}^{9} B_{k} \lambda^{9-k}, \\
% %
% [{W^6}_{\la}W^6] &= \sum_{k=0}^{11} C_{k} \lambda^{11-k}.
% \end{split}
% \end{align}
\begin{align}\label{eq:generalOPE.case1}
\begin{split}
[{W^4}_{\la}W^4] = \sum_{k=0}^{7} A_{k} \lambda^{7-k}, \qquad
[{W^4}_{\la}W^6] = \sum_{k=0}^{9} B_{k} \lambda^{9-k}, \qquad
[{W^6}_{\la}W^6] = \sum_{k=0}^{11} C_{k} \lambda^{11-k}.
\end{split}
\end{align}
Theorem \ref{thm:PBW} asserts that $W^k(B_3)$ has a basis consisting of PBW monomials in generators $L$, $W^4$ and $W^6$ and their derivatives. We would like to know that $W^k(B_3)$ is unique with this property. For this we carry out the basic strategy on \eqref{eq:generalOPE.case1}, equating coefficients in \eqref{eq:comm.fla}, and solving (now there are $216$ variables instead of the $75$ of the previous two cases). Indeed this system has a unique solution.

In relation to the previous two cases, there is now an additional complication:  the PBW monomials in $L$, $W^4$ and $W^6$ of conformal weight at most $11$ are not linearly independent in the vertex algebra $W$. Indeed by comparing characters we see that the kernel of the quotient map ${W^k(B_3)} \rightarrow {W(B_3, 8/7)}$ is trivial in conformal weight up to $7$ and has a one-dimensional component in conformal weight $8$. We confirm that the OPE relations of $W^k(B_3)$ obtained above, do indeed admit a singular vector of the form $W^4_{(-1)}W^4 - \sigma_8$ where
\begin{align}\label{eq:W4W4=}
\sigma_8 \in V(c)_8 \oplus M(c, 4)_8 \oplus M(c, 6)_8.
\end{align}

We must now check that for triples of $L$, $W^4$ and $W^6$ in the quotient $W$, the equation \eqref{eq:comm.fla} admits a unique solution. This does not follow from the uniqueness results of the previous paragraphs. Indeed as remarked in {\cite[Section 3]{Linshaw-W-infty}}, not every finitely strongly generated vertex algebra is a quotient of one with a PBW basis (the commutant of the Heisenberg subalgebra of $V^k(\sll_2)$ being a prominent example).

Keeping this subtlety in mind, the basic strategy nevertheless remains the same as in cases (2) and (3) considered above: expand \eqref{eq:comm.fla} for all triples taken from $L$, $W^4$ and $W^6$ and equate coefficients relative to a monomial basis, adding now an equation $W^4_{(-1)}W^4 - \sigma_8 = 0$ where $\sigma_8$ of the form \eqref{eq:W4W4=} is taken as an unknown linear combination of monomial Virasoro descendents of $\vac$, $W^4$ and $W^6$ (this requires $14$ additional coefficients). We retain only those equations coming from conformal weight at most $8$. Fortunately, up to normalisation of $W^4$ and $W^6$ by scalars, the equations are sufficient to constrain the coefficients in \eqref{eq:generalOPE.case1} and \eqref{eq:W4W4=} to a unique solution.

From all this we deduce that there exists a homomorphism of vertex algebras
\begin{align*}
\varphi : W^k(B_3) \rightarrow W.
\end{align*}
Since we have computed $\chi_{W}(q) \equiv \chi_{W_k(B_3)}(q) \bmod{q^{25}}$, we conclude as in other cases that
\[
W = \im(\varphi) \cong W_k(B_3) = W(B_3, 8/7).
\]
%just as in the other cases.

\begin{rem}
In \cite{K.Sun} K. Sun has conjectured an isomorphism
\[
W(E_7(a_1), 19/16) \cong \left(\vir^{N=1}_{16, 2}\right)_{\ov 0},
\]
where the right hand side is the even part of the $(16, 2)$ minimal model over the $N=1$ superconformal algebra.
\end{rem}

\end{document}